\documentclass[smallcondensed,nospthms]{svjour3}

\usepackage[utf8]{inputenc}
\def\final{1}
\usepackage{color}
\usepackage{subfigure}
\usepackage{url}
\usepackage{graphicx}
\usepackage{amsmath,amssymb}
\usepackage{amsthm}
\usepackage{bm}
\usepackage{pinlabel}
\usepackage[plainpages=false]{hyperref}

\usepackage{mathptmx}
\usepackage[square,numbers,sort&compress]{natbib}

\ifnum\final=0
\newcommand{\mynote}[1]{\marginpar{\tiny\sf #1}}
\else
\newcommand{\mynote}[1]{}
\fi

\usepackage{theomac}
\newtheoremWithMacro{theorem}{Theorem}[section]

\newtheorem{lemma}{Lemma}[section]

\newtheorem{question}{Question}

\newtheoremWithMacro{corollary}{Corollary}
\newtheoremWithMacro{prop}[lemma]{Proposition}

\newcommand{\figref}[1]{Figure \ref{fig:#1}}

\newcommand{\lemref}[1]{Lemma \ref{lemma:#1}}

\newcommand{\propref}[1]{Proposition \ref{prop:#1}}
\newcommand{\theoref}[1]{Theorem \ref{theo:#1}}
\newcommand{\theorefX}[1]{\ref{theo:#1}}
\newcommand{\corref}[1]{Corollary \ref{cor:#1}}
\newcommand{\secref}[1]{Section \ref{sec:#1}}

\newcommand{\lemlab}[1]{\label{lemma:#1}}
\newcommand{\proplab}[1]{\label{prop:#1}}
\newcommand{\theolab}[1]{\label{theo:#1}}
\newcommand{\corlab}[1]{\label{cor:#1}}
\newcommand{\seclab}[1]{\label{sec:#1}}

\newcommand{\R}{\mathbb{R}}
\newcommand{\Z}{\mathbb{Z}}

\newcommand{\HH}{\mathrm{H}}

\newcommand{\eop}{\hfill $\qed$}

\begin{document}
\title{Topological designs}
\author{Justin Malestein 	\and	Igor Rivin 	\and	Louis Theran
}
\institute{
Justin Malestein \at
Math Department, Hebrew University, Jerusalem, Israel \\
\email{justinmalestein@math.huji.ac.il}
\and
Igor Rivin \at
Math Department, Temple University, Philadelphia, PA 19122 \\
\email{rivin@math.temple.edu}
\and
Louis Theran \at
Institut für Mathematik, Freie Universität Berlin, 14195 Berlin, Germany \\
\email{theran@math.fu-berlin.de}
}
\date{}
\maketitle
\begin{abstract}
We give an exponential upper and a quadratic lower bound on the number of pairwise non-isotopic simple closed
curves which can be placed on a closed surface of genus $g$ such that any two of the curves intersects at most once.
Although the gap is large, both bounds are the best known for large genus.  In genus one and two, we
solve the problem exactly.

Our methods generalize to variants in which the allowed number of pairwise
intersections is odd, even, or bounded, and to surfaces with boundary components.
\end{abstract}

\section{Introduction} \seclab{intro}
The classical area of \emph{combinatorial designs} (see, e.g., \cite{S04}) asks about
extremal subset systems of a finite set that have specified intersection patterns.
In this paper, we are interested in extremal families of curves in a closed, orientable
surface with specified intersection patterns.  More specifically, Benson Farb and Chris Leininger
brought the following topological variant to our attention:
\begin{question} \label{farbq}
Let $S$ be a closed oriented surface of genus $g$, and let $X$ be a collection of pairwise
non-isotopic essential simple closed curves $\gamma_1, \dotsc, \gamma_N$ such that $\gamma_i$ intersects
$\gamma_j$ at most once for any $i\neq j$. How large can $X$ be?
\end{question}
In fact, this question goes back further, to work of Juvan, Malni{\v{c}}, and Mohar \cite{JMM96}.
Juvan, et al. define a \emph{$k$-system} of curves on a genus $g$ closed, oriented surface to be a
collection of non-isotopic, essential, simple closed curves $\gamma_1, \ldots, \gamma_N$ such that each of the
$\gamma_i$ intersects any $\gamma_j$, $j\neq i$ at most $k$ times. We define $N(k,g)$ to be the maximum
size of any $k$-system on a closed surface $S$ of genus $g$.  The most important related results from
\cite{JMM96} are:
\begin{itemize}
\item For all $g\ge 0$ and $k\ge 0$, $N(g,k)$ is finite.
\item $N(k,1)\le 2k + 3$.
\item $N(k,g)\ge (k/g)^g$
\end{itemize}

\subsection{Results and techniques}
We study $N(k,g)$ for $k$ fixed.  We are particularly interested in the
Farb-Leininger case $k=1$, for which we can prove:
\begin{theorem}\theolab{generalbounds}
For $g=3$,
\[
6g \le N(1,g) \le (g-1)(2^{2g} -1)
\]
and for all $g\ge 4$,
\[
g^2+\frac{5}{2}g\le N(1,g)\le (g-1)(2^{2g} -1)
\]
\end{theorem}
Both the upper and lower bounds are, to the best of our knowledge, the best known.  Juvan, et al.
don't optimize the bound in \cite[Theorem 3.3]{JMM96} for $k=1$, but examining their arguments, we get an
upper bound on the order of $g!$.

The upper bound of \theoref{generalbounds} is based on the $\Z/2\Z$-homology
of curves (\propref{mod2classupper}).
Thus, the same upper bound statement and proof applies if we replace ``$1$-system'' with
$k$-system in which all pairs of curves are either disjoint or intersect an odd number of times.
Thus the upper bound in \theoref{generalbounds} has some built-in slack, and it is plausible that,
for large enough $k$, this upper bound is nearly tight.  For large enough $k$ and $g$,
\cite[Proposition 5.1]{JMM96} gives a lower bound for $N(k,g)$ on the order of $2^{c\sqrt{g}}$.
Reducing the gap will require a proof method that can see the difference between
``one intersection'' and ``an odd number of intersections''.

\paragraph{Surfaces with boundary components}
For a surface with boundary components, a $k$-system is defined similarly to the closed case, except
we add the condition that none of the curves is allowed to to be homotopic to a boundary curve.
Define $N(k,g,n)$ to be the maximum size of a $k$-system in a surface of genus $g$ with $n$ boundary components.

Holding $g$ fixed (along with $k$), we can determine $N(1,g,n)$ very precisely in terms of $N(1,g)$:
\begin{theorem} \theolab{boundaryasymptotics}
Fix some $g \geq 2$.  Then, for all $n \geq 0$,
$$ (2g+1)(n+1) \leq N(1, g, n) \leq N(1,g) + (2g+1)n$$
For $g = 1$ and $n \geq 1$
$$N(1, 1, n) =3 \cdot n $$
\end{theorem}
Relaxing the definition of a $k$-system to allowed arcs as well as curves on a surface with
boundary seems to make the problem as hard as determining $N(1,g)$.

\paragraph{Low genus}
On the torus, it is straightforward to show that a maximal $1$-system has $3$ curves, i.e. $N(1,1) = 3$.  In genus two, we present
an elegant geometric argument, based on the hyperelliptic involution, which shows that:
\begin{theorem}\theolab{genus2}
N(1,2) = 12.
Moreover, there are exactly two mapping class orbits of maximal $1$-systems on a genus-$2$ surface.
\end{theorem}
The proof gives exact structural information about the maximal $1$-systems on a genus-$2$ surface:
they are lifts of $6$-vertex triangulations of the sphere drawn in the quotient of the surface by the
hyperelliptic involution.

\paragraph{Nearly all curves intersecting}
Given a $k$-system $X=\{\gamma_{1},\gamma_2,\ldots, \gamma_N\}$ on a closed surface $S$ of genus $g$, we define its
\emph{odd intersection graph} $G_{\text{odd}}(X)$ to be the graph that has one vertex for each of the $\gamma_i$ and an edge
between vertices $i$ and $j$ if $\gamma_i$ intersects $\gamma_j$ an odd number of times.  With assumptions on
$G(X)$, we can get much tighter upper bounds.
\begin{theorem}\theolab{alloddupper}
Let $X=\{\gamma_{1},\gamma_2,\ldots, \gamma_N\}$ be a $k$-system with odd intersection graph $G_{\text{odd}}(X)$.  Then for $g\ge 1$
and $k\ge 1$:
\begin{itemize}
\item If $G_{\text{odd}}(X)$ is the complete graph $K_N$, then $N\le 2g+1$, and the bound is sharp for all $k \geq 1$.
\item If $G_{\text{odd}}(X)$ has average degree $N-D$, then $N\le (D+1)(2g+1)$.
\end{itemize}
\end{theorem}
The first statement is proved using a $\Z/2\Z$-homology argument and the second statement follows via an application of
Turán's Theorem \cite{T41} to $G(X)$'s complement, an idea we learned from Van Vu's paper \cite{V99}.

\paragraph{Nearly all curves non-intersecting}
Given a $k$-system $X=\{\gamma_{1},\gamma_2,\ldots, \gamma_N\}$ we define its \emph{intersection graph}
to be the graph that has one vertex for each $\gamma_i$ and an edge between vertices $i$ and $j$ if $\gamma_i$
intersects $\gamma_j$.

\begin{theorem}\theolab{lowdegree}
Let $X=\{\gamma_{1},\gamma_2,\ldots, \gamma_N\}$ be a $k$-system on a closed surface of genus $g$
with  intersection graph $G(X)$. For $g\ge 2$ and $k\ge 1$.  If $G(X)$ has average degree $D$, then $|X|\le (D+1)(3g-3)$.
\end{theorem}
In the case where $D = 0$, this is just the well-known fact that a maximal system of disjoint nonisotopic simple closed
curves has at most $3g-3$ curves.

\subsection{Acknowledgements}
JM, IR, and LT received support for this work from Rivin's NSF CDI-I grant DMR 0835586.
LT's final preparation was supported by the European Research Council under the European Union’s
Seventh Framework Programme (FP7/2007-2013) / ERC grant agreement no 247029-SDModels.  Our
initial work on this problem took place at Ileana Streinu's 2010 Barbados workshop at
the Bellairs Research Institute.

\section{Genus $2$} \seclab{gen2}
In this section, we give a geometric argument to prove \theoref{genus2}.
The argument generalizes to give a sub-quadratic lower bound in higher genus.

\subsection{The hyperelliptic involution}
The key tool we need is the \emph{hyperelliptic involution}.  Here are the facts
we need, which can all be found in \cite{HS89}.
\begin{prop}
\proplab{hyper-exists}
Every hyperbolic surface $S$ of genus $2$	admits a ``hyperelliptic involution'', acting by an isometry,
with $6$ fixed points.
\end{prop}
The fixed points of the hyperelliptic involution are called \emph{Weierstrass points}, and
we denote them $w_1, w_2, \ldots, w_6$.
\begin{prop}
\proplab{hyper-quotient}
Let $S$ be a closed hyperbolic surface of genus $2$.  The quotient $Q$  of $S$ by the hyperelliptic
involution is a sphere with six double points corresponding to the Weierstrass points; i.e., $Q$
is an orbifold of signature $(0; 2, 2, 2, 2, 2, 2)$.
\end{prop}
Simple closed geodesics in $S$ have nice representatives in $Q$.
\begin{prop}
\proplab{hyper-geodesics}
Let $S$ and $Q$ be as in \propref{hyper-quotient}.  Then any non-separating simple closed geodesic $\gamma$ in
$S$ goes through exactly two of the Weierstrass points and projects to a geodesic segment going between
exactly two of the doubled points in $Q$.
\end{prop}

\subsection{Proof of \theoref{genus2}}
By \propref{nosep} proven below, there is no separating curve in a maximal $1$-system, and by
\propref{hyper-geodesics}, the $1$-systems in $S$ that have only non-separating curves
correspond to isomorphism classes of simple planar graphs on $6$ vertices.
Since planar triangulations on $6$ vertices have $12$ edges, it follows that $N(1,2) = 12$.

There are two graph isomorphism classes of planar triangulations%
\footnote{This can be checked using Brendan McKay's {\tt plantri} \cite{BM07}.}
on $6$ vertices, namely, the octahedron and doubly stellated tetrahedron.
Since these graphs are $3$-connected, there are only two triangulations
up to homeomorphism of the sphere (permuting Weierstrass points).
Birman and Hilden \cite{BH71} have shown that the mapping class group
of the genus $2$ surface modulo the hyperelliptic involution is naturally isomorphic to
the mapping class group of the $6$-times punctured sphere.
Consequently, there are two mapping class group orbits of maximal $1$-systems.
\eop

\subsection{A linear lower bound}
The lower bound of \theoref{genus2} generalizes to surfaces $S$ of $g\ge 3$: there are $2g+2$
Weierstrass points corresponding to double points in the quotient by the hyperelliptic involution,
implying that any planar graph on $2g+2$ vertices corresponds to some $1$-system.  Thus, we
obtain the lower bound for $g = 3$ from \theoref{generalbounds}.

\section{Structure of maximal $k$-systems} \seclab{structure}
We will need the following result on the structure of maximal $1$-systems:
\begin{prop} \proplab{nosep}
Let $S$ be a closed surface of genus $g\ge 1$.  Then, any maximal $1$-system in $S$ contains
no separating curves.
\end{prop}
From now on in this section, we assume that $S$ satisfies the hypotheses of \propref{nosep}.
We also use the standard notions of \emph{minimal position} and \emph{geometric intersection number}
for curves, which may be found in \cite[Section 1.2.3]{FM12}.  We denote the geometric
intersection number of (the isotopy classes of) $\gamma_1$ and $\gamma_2$ by $i(\gamma_1,\gamma_2)$.
A $k$-system is in minimal position when the curves in it are pairwise in minimal position.  Every
$k$-system has such a representative \cite[Corollary 1.9]{FM12}.

The auxiliary \lemref{annulipoly}, which describes the complementary regions of a maximal $k$-system is also
of independent interest. Let $X$ be a $k$-system and define the \emph{complementary regions of $X$} to be
the components of the surface $S_c$ obtained by cutting $S$ along the curves in $X$.
\begin{lemma} \lemlab{annulipoly}
Let $X$ be a maximal $k$-system such that the curves in $X$ are in minimal position.
Then, the complementary regions are polygons and annuli. Furthermore, each annulus has a
boundary component consisting of a single curve in $X$.
\end{lemma}
Intuitively, one would expect that if \propref{nosep} were to fail, an inductive argument would
establish an $O(g)$ upper bound on the size of any $1$-system, contradicting the fact that
there are $1$-systems with $\Omega(g^2)$ curves.  However, we are unaware of a rigorous proof along
these lines, so instead we use \lemref{annulipoly}.

The difficult part of the proof of \lemref{annulipoly} will be to rule out pairs of pants
as complementary regions.  This next lemma is straightforward and gives the
starting point.
\begin{lemma}\lemlab{nobiggerthanpants}  Suppose $X$ is a maximal
$k$-system in minimal position.  Then, the complementary regions are all polygons,
annuli, or pairs of pants.
Furthermore, each annulus has a boundary component consisting of a single curve in
$X$, and every boundary component of the pairs of pants are curves in $X$.
\end{lemma}
\begin{proof}
If a complementary region has nonzero genus, then it contains a simple closed
curve that is not isotopic to a boundary curve which we could add to $X$ to obtain a
larger $k$-system.  A similar fact is true if a complementary region is a sphere with at least
$4$ boundary components.  This proves that the complementary regions are all polygons, annuli, or
pairs of pants.

Now consider a complementary region that is an annulus.  If neither boundary component is a single
curve in $X$, then the core curve of this annulus is, by our assumption of minimal position, not
isotopic to any curve in $X$ and does not intersect any curves in $X$.  This contradicts the
maximality of $X$.

Similar reasoning shows that all boundary curves of a pair of pants must be elements of $X$.
\end{proof}

\subsection{Proof of \lemref{annulipoly}}
Continuing from \lemref{nobiggerthanpants}, it suffices to rule out pants as a complementary region.
Comparing the (classical) upper bound, \theoref{all-disjoint}, of the size of a $0$-system and the
lower bound on the size of a $1$-system from \theoref{generalbounds},
we see that not all the complementary regions can be pairs of pants.

Thus, if $X$ has some pair of pants as a complementary region, then
connectedness of $S$ and \lemref{nobiggerthanpants} implies that some curve $\gamma$ in $X$ is disjoint
from the rest of $X$ and bounds a pair of pants $P$ on one side and an annulus $A$ on the other.  The other boundary
component of $A$ consists of arcs of curves in $X$, and there are at least two such arcs.  Let $\alpha, \beta$ be
two curves with consecutive arcs in $\partial A$.
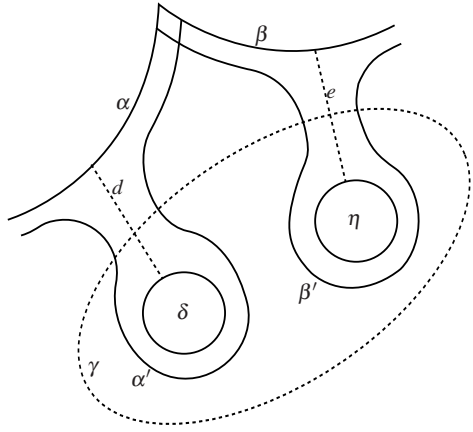
\begin{figure}[htbp]
\begin{center}
\labellist \small \hair 2pt

\pinlabel $\alpha$ at 50 140
\pinlabel $\alpha'$ at 58 21
\pinlabel  $\beta$ at 110 170
\pinlabel $\beta'$ at 130 58
\pinlabel $\gamma$ at 37 25
\pinlabel $\delta$ at 76 50
\pinlabel $\eta$ at 150 89
\pinlabel $d$ at 48 104
\pinlabel $e$ at 140 145
\endlabellist

\includegraphics[width=0.5\textwidth]{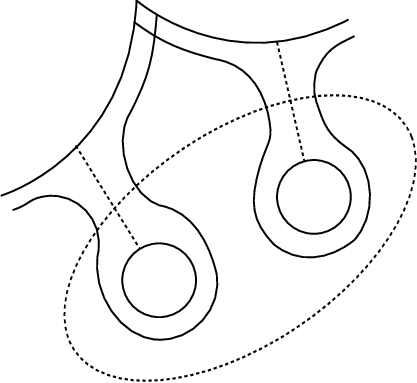}
\end{center}
\caption{A depiction of curves $\alpha, \alpha', \beta, \beta', \gamma, \delta, \eta$ and arcs $d, e$.}
\label{figure:newcurves}
\end{figure}

We now show that $\gamma$ may be replaced with two new curves to obtain a larger $1$-system; the construction is
depicted in Figure \ref{figure:newcurves}.  Let $\delta, \eta$ be
the other boundary components (aside from $\gamma$) in the pair of pants $P$, and let $d$ be a simple arc connecting
$\alpha$ to $\delta$ and lying in $A \cup P$.
The regular neighborhood of $d \cup \alpha \cup \delta$
is a pair of pants where two of the boundary curves are isotopic to $\alpha$ and $\delta$; let $\alpha'$ be the other
boundary component.  Similarly, let $e$ be a simple arc which connects $\beta$ to $\eta$, lies in $A \cup P$,
and is disjoint from $d$;  let $\beta'$ be defined similarly to $\alpha'$.

To show $Y = X \cup \{\alpha', \beta' \} \setminus \{\gamma\}$ is a $k$-system of larger size, we must show pairwise intersections are less than or equal to $k$ and that
$\alpha'$ and $\beta'$ are not isotopic to other curves in $X$ or each other.  It is possible to do the latter by establishing the following geometric intersection numbers
since the numbers are an invariant of isotopy classes.

\begin{itemize}
\item[(1)] For all curves $\mu \in X \setminus \{\gamma\}$, we have $i(\alpha', \mu) = i(\alpha, \mu)$ and $i(\beta', \mu) = i(\beta, \mu)$
\item[(2)] $i(\alpha', \beta') \leq k$
\item[(3)] $i(\alpha', \gamma) = 2 = i(\beta', \gamma)$
\end{itemize}

Before proving (1)-(3), let us see that they establish that $Y$ is a $k$-system.  Statements (1) and (2) show curves in $Y$ intersect at most $k$ times.  (Note that (1) also says
$i(\alpha', \alpha) = i(\alpha, \alpha) = 0$ and similarly for $\beta, \beta'$.)  Statement (3) implies $\alpha'$ and $\beta'$ are not isotopic to any curves
in $X$ since they are all disjoint from $\gamma$, and statement (1) implies $i(\alpha', \alpha) = 0 \neq i(\alpha, \beta) = i(\alpha, \beta')$ which
establishes that $\alpha'$ and $\beta'$ are not isotopic.

\paragraph{Proof of (1)} We only show that $i(\alpha', \mu) = i(\alpha, \mu)$ since the proof applies also to $\beta$ and $\beta'$.  If $\mu = \alpha$,
this is immediate since $\alpha'$ and $\alpha$ are disjoint.  Now suppose $\mu \neq \alpha$.

\begin{figure}[htbp]
\begin{center}
\includegraphics[width=0.5\textwidth]{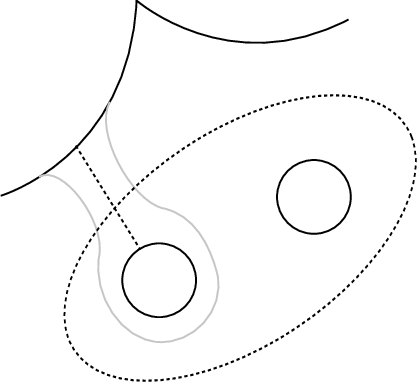}
\end{center}
\caption{The arc $a''$ indicated in gray.}
\label{figure:nopantsfig2}
\end{figure}

Homotope $\alpha'$ to a simple closed curve $\alpha''$
which is the union of an arc $a$ from $\alpha$ and an arc $a''$ which lies in $A \cup P$, is disjoint from all curves in $X \setminus \{\gamma\}$,
and cuts $P$ into two annuli.  (See Figure
\ref{figure:nopantsfig2}.)  Clearly, the arc
$a$ and curve $\mu$ cross $i(\alpha, \mu)$ times, so by the bigon criterion (\cite[Proposition 1.7]{FM12}), it suffices to show that it is impossible to bound a disc
with an arc from $\mu$ and an arc from $\alpha''$.  Since $X$ is in minimal position, $\alpha$ and $\mu$ have no bigons, and so any bigon made from
$\mu$ and $\alpha''$ must use the arc $a''$.  However, on each side of $a''$, we may connect $a''$ to either $\delta$ or $\delta'$ with an arc disjoint from $\mu$
(since $\mu$ is disjoint from $P$);
this would be impossible if an arc containing $a''$ made a bigon with an arc from $\mu$.  This established the equality.

\paragraph{Proof of (2)} Let $\alpha''$ be as above.  By further homotoping, we can ensure that $a''$ is disjoint from $\beta'$ and so
$$i(\alpha', \beta') = i(\alpha'', \beta') \leq i(\alpha, \beta') = i(\alpha, \beta) \leq k$$

\paragraph{Proof of (3)}
Again, we use the bigon criterion.  By construction, $\gamma$ and $\alpha'$ each cut each other into two arcs, $c_1, c_2$ and $a_1', a_2'$ respectively, where
w.l.o.g. $a_1'$ lies in $P$.  The curves $a_1' \cup c_1$ and $a_1' \cup c_2$ are each homotopic to some boundary component in $P$ and hence cannot be trivial.
If $a_2' \cup c_i$ bounded a disc, then we could homotope $\alpha'$ to lie in $P$ which is impossible since $\beta$ and $\alpha'$ have nonzero geometric intersection
number.
\eop
\subsection{Proof of \propref{nosep}}
Suppose, for a contradiction, that $X$ is a maximal $1$-system with a separating curve $\gamma$;
w.l.o.g., we may assume that $X$ is in minimal position.  \lemref{annulipoly} then applies,
so $\gamma$ is incident on two complementary regions $A$ and $A'$ that are both annuli,
and each of $A$ and $A'$ has a boundary component consisting of arcs from at least two different
curves in $X$.  We will show that by removing $\gamma$, two new curves may be
added to get a larger $1$-system.

Because $X$ is a $1$-system, we may slightly strengthen the conclusion of \lemref{annulipoly}:
the boundary components of $A$ and $A'$ that are not $\gamma$ span arcs from at least
three different curves in $X$.  Let $\alpha, \beta$ be curves
with consecutive arcs in $\partial A$, and let $\delta$ be a curve contributing a boundary arc to $A'$. Let $S'$ be the surface
obtained by cutting along $\gamma$. There is an arc $a$ starting and ending at $\partial S'$ such that $a$ ``follows'' $\alpha$. More
precisely, there is an arc $a$ such that $a$, $\alpha$ and a part of $\partial S'$ cobound an annulus and $a$ has the same intersection
numbers with the other curves in $X \setminus \{\gamma\}$. Similarly, there is such an arc $b$ for $\beta$, and let $d_1, d_2$ be two
disjoint such arcs for $\delta$. We can homotope the arcs $a, b$ so that:%
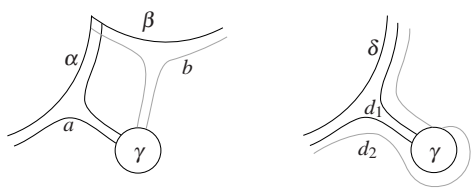
\begin{figure}[htbp]
\begin{center}
\labellist \small \hair 2pt

\pinlabel $\alpha$ at 50 100
\pinlabel  $\beta$ at 110 130
\pinlabel $\gamma$ at 102 26
\pinlabel $\gamma$ at 331 26
\pinlabel $\delta$ at 284 110
\pinlabel $a$ at 48 48
\pinlabel $b$ at 140 95
\pinlabel $d_1$ at 285 62
\pinlabel $d_2$ at 280 30
\endlabellist
\includegraphics[width=0.5\textwidth]{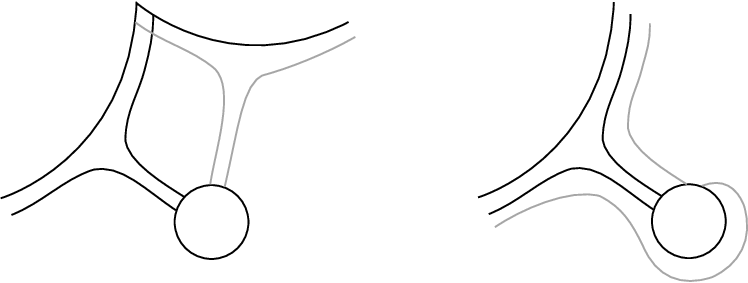}
\end{center}
\caption{A depiction of the arcs to be connected on either side of the separating curve.
Note that the gray arcs will glue together.}
\label{figure:arcsfornosep}
\end{figure}
\begin{itemize}
\item $a$ and $b$ intersect exactly once transversely and
\item the endpoints of $a$ match up with $d_1$ and $b$ with $d_2$ when gluing $S'$
back together.
\end{itemize}
(This is depicted in Figure \ref{figure:arcsfornosep}.)

Let $\alpha'$ and $\beta'$ be the resulting arcs from gluing $a$ with $d_1$ and $b$ with $d_2$ respectively.
Since $\alpha'$ and $\beta'$ intersect exactly once transversely, they are distinct non-trivial simple closed curves.
Furthermore both $\alpha'$ and $\beta'$ intersect some curve exactly once on either side of $\gamma$. Since no other
curve in $X$ does this, $\alpha'$ and $\beta'$ are isotopically distinct from all other curves in $X$. Thus
$X' = X \cup {\alpha', \beta'} \setminus \{\gamma\}$ is a larger $1$-system than $X$.
\eop

\section{Upper bounds from homology} \seclab{upper-hom}
Now we turn to the general case of a surface $S$ of genus $g\ge 2$.  The point of this section is
to prove the following proposition, from which the upper bound in \theoref{generalbounds}
follows readily.

\begin{prop}\proplab{mod2classupper}
Let $S$ be a closed oriented surface of genus $g\ge 2$.  Then any $1$-system of curves on $S$
has at most $g-1$ curves in any nontrivial $\Z/2\Z$-homology class.
\end{prop}

\subsection{All curves null-homologous mod $2$}
The main lemma we need is:
\begin{lemma}\lemlab{disj-all-sep}
Let $g\ge 1$ and let $S$ be a genus $g$ surface with $2$ boundary components.  Let $X$ be a $0$-system  of separating curves
in $S$ such that each curve separates the boundary components.  Then, $|X| \le g-1$.
\end{lemma}
\begin{proof}
We argue by induction.  Suppose $g = 1$.  In this case, the only separating simple curves which are not boundary parallel
cut $S$ into a torus with a single boundary component and a three-holed sphere.  Such a curve does not separate the
boundary components of $S$.

Suppose $g > 1$.  Cutting $S$ along a curve $\gamma$ in $X$ yields
two surfaces $S', S''$ with genus $g', g''$ adding up to $g$ and each with $2$ boundary components.
We note that $X \setminus \{\gamma\}$ deposits a $0$-system in each of $S'$ and
$S''$ with the same properties as $X$.  Thus $|X| \le 1 + g' - 1+ g''-1 = g - 1$.
\end{proof}

\subsection{Proof of \propref{mod2classupper}}
It suffices to show that one can find at most $g-1$ mutually disjoint curves in the same $\Z/2\Z$-homology class.
Let $S$ be a closed genus $g\ge 2$ surface, and $X$ be a $0$-system of curves all in the same
$\Z/2\Z$-homology class.

We reduce the proposition to \lemref{disj-all-sep}.  Cut along some $\gamma\in X$ to obtain a surface $S'$ of genus
$g' = g -  1 \geq 1$ with two boundary components corresponding to $\gamma$;  fill in each boundary component with a
disk to get a closed surface $S''$. The curves in $X\setminus \{\gamma\}$ in $S''$ are all null-homologous in
$\HH_1(S'',\Z/2\Z)$.  This means that their homology classes in $\HH_1(S'',\Z)$ are all either non-primitive or trivial.  The
former implies the the classes are not primitive in $\HH_1(S,\Z)$, which is disallowed by the hypothesis that all
the curves in $X$ are simple (see \cite[Proposition 6.2]{FM12}).
Thus, any curve in $X \setminus \{\gamma\}$ separates $S''$.

Furthermore, any $\eta \in X \setminus \{\gamma\}$ must separate the boundary components of $S'$.  If it did not,
then it would bound a subsurface in $S'$ with a single boundary component,
and thus also in $S$.  This however would imply that $\eta$ is null-homologous as an element of $\HH_1(S, \Z/2\Z)$,
a contradiction.  Consequently, we are in the situation of \lemref{disj-all-sep} and are done.
\eop

\subsection{Proof of the upper bound from \theoref{generalbounds}}
By \propref{nosep}, a maximal $1$-system contains no separating curve and thus no null-homologous curve.  The curves must then all lie in the $2^{2g} -1$ nontrivial
$\Z/2\Z$-homology classes, and, by \propref{mod2classupper},
each of these contains at most $g-1$ curves in a $1$-system.
\eop

\section{(Nearly) all intersecting or disjoint} \seclab{extreme}
In this section we give much sharper upper bounds when the intersection graph $G(X)$ of a
$1$-system $X$ is either very sparse or very dense.  The key cases are when $G(X)$
or its complement are complete.

\subsection{All curves disjoint}
We recall the following classical and widely known fact, which may be found in
\cite[Section 8.3.1]{FM12}.
\begin{theorem}
\theolab{all-disjoint}
Let $g>1,$ and let $S$ be a closed surface of genus $g$, and let $X$ be a $0$-system.
Then $|X|\leq 3g - 3$.
\end{theorem}

\subsection{All curves intersecting}
To bound the size of $1$-systems with all pairs of curves intersecting, we continue
along the lines of \secref{upper-hom}.  Let $v$ and $w$ be vectors in $(\Z/2\Z)^{2g}$
and let $(v,w)$ denote the standard symplectic pairing.
\begin{prop}\proplab{linalg}
Let $v_1,v_2,\ldots, v_N$ be non-zero vectors in $(\Z/2\Z)^{2g}$ with the property that, for all $i\neq j$,
$(v_i,v_j) = 1$.  Then $N\le 2g+1$.
\end{prop}
\begin{proof}
Suppose there is a linear dependence
\begin{equation}
\label{lin-dep}
\alpha_1 v_1 + \alpha_2 v_2 + \cdots + \alpha_k v_k = 0
\end{equation}
among the $v_i$, with not all the $\alpha_i$ zero. Suppose further that the dependence \eqref{lin-dep}
is non-trivial and that $\alpha_1=1$ and $\alpha_2=0$.
Pairing both sides of \eqref{lin-dep} with $v_1$ tells us the number of non-zero $\alpha_i$
is odd; similarly pairing both sides of \eqref{lin-dep} with $v_2$ tells us the number of non-zero
$\alpha_i$ is even.  The resulting contradiction implies that, in fact, for any non-trivial linear
dependence among the $v_i$, the $\alpha_i$ are all one.

Thus, the $v_i$ are either independent, implying $N\le 2g$,
or there is a unique linear dependence with full support among them,
implying $k=2g+1$.
\end{proof}
As a corollary, we obtain:
\begin{prop}\proplab{all-intersecting}
Let $k\ge 1$, and let $X$ be a $k$-system in a closed genus $g$ surface with any pair of
curves intersecting an odd number of times.  Then $|X|\le 2g + 1$.
\end{prop}
\begin{proof}
If the minimal geometric intersection number between essential, non-isotopic, simple closed
curves is odd, then their algebraic intersection number mod $2$ is $1$ and, furthermore,
their $\Z/2\Z$-homology classes must be distinct.  Apply \propref{linalg}.
\end{proof}

\subsection{Proofs of Theorems \theorefX{alloddupper}--\theorefX{lowdegree}}
We prove \theoref{alloddupper}, since the proof of \theoref{lowdegree}
is nearly identical.  The upper bound when the graph $G_{\text{odd}}(X)$
is $K_N$ is \propref{all-intersecting}.  If $G_{\text{odd}}(X)$ has
average degree $N-D$, then its complement has average degree $D$, and
by Turán's Theorem \cite{T41}, must contain an independent set
of size $N/(D+1)$.  Applying \propref{all-intersecting} to the curves in $X$ represented by the
corresponding clique in $G_{\text{odd}}(X)$, we see that $2g + 1\ge N/(D+1)$.
\eop

\subsection{\propref{linalg} is sharp}
The bound of \propref{linalg} is tight, as the following example shows.  Define
$v^1_1  =  (1,1)$, $v^1_2  =  (0,1) $, and $v^1_3  = (1,0)$.  For $g\ge 2$,
inductively define $v^g_{i}  =  (v^{g-1}_{i} ; 0, 0)$ for $i\in [1,2g-1]$ and
then $v^g_{2g -1}  = (v^{g-1}_{2g - 1} ; 1, 0)$, $v^g_{2g}  = (v^{g-1}_{2g - 1} ; 0, 1)$,
and $v^g_{2g}  = (v^{g-1}_{2g + 1} ; 1, 1)$.  (The semi-colons mean concatenating vectors.)
\begin{figure}[htbp]
\centering
\includegraphics[width=0.7\textwidth]{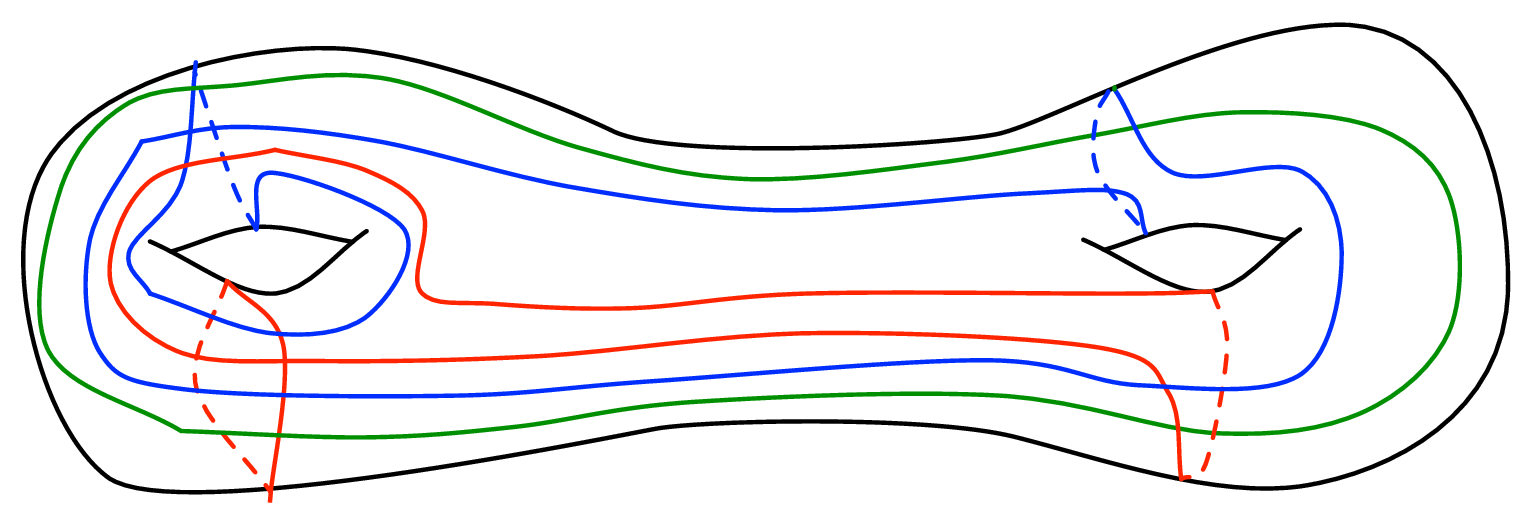}
\caption{Curves corresponding to $v^2_i$.}
\label{fig:curves-g2}
\end{figure}
\begin{prop}\proplab{linalg-tight}
Let $g\ge 1$, and let $v^g_i$ be defined as above.  For all $j\neq i$, we have
$(v^g_i,v^g_j) = 1$.
\end{prop}
\begin{proof}
For $g=1$, this is an easy computation, and the $g\ge 2$ cases
follow by induction.
\end{proof}
\figref{curves-g2} shows curves with $\Z/2\Z$-homology classes given
by the $v^2_i$. The vectors $v^g_i$ are canonical.

\begin{prop}\proplab{linalg-canonical}
Let $g\ge 1$ and suppose that $w_1, w_2, \ldots, w_{2g+1}$ have the property that,
for $i\neq j$, $(w_i,w_j)=1$.  Then there is a symplectic automorphism $A$ of $(\Z/2\Z)^{2g}$
such that $Av_i = w_i$.
\end{prop}
In the proof, we need the basic fact:
\begin{lemma}\lemlab{linalg-indep}
Let $v_1,v_2,\ldots,v_{2g}$ and $w_1,w_2,\ldots,w_{2g}$ be linearly independent
sets of vectors such that $(v_i,v_j) = (w_i,w_j)$ for all $i$ and $j$.
Then there is a symplectic automorphism $A$ of $(\Z/2\Z)^{2g}$
such that $Av_i = w_i$.
\end{lemma}

\begin{proof}[Proof of \propref{linalg-canonical}]
Let $A$ be the automorphism from \lemref{linalg-indep} applied to
the $v^g_1,v^g_2,\ldots, v^g_{2g}$ and $w_1,\ldots,w_{2g}$ from the statement of the Theorem.
This is allowed, because the proof of \propref{linalg} says that for any collections
meeting the hypothesis of the $v^g_i$ and $w_i$, we have
\begin{eqnarray*}
v^g_{2g+1} & = & v^g_1 + v^g_2 + \cdots + v^g_{2g} \\
w_{2g+1} & = & w_1 + w_2 + \cdots + w_{2g}
\end{eqnarray*}
with the first $2g$ vectors independent.  Thus, we see that
\begin{equation*}
Av^g_{2g+1} = A( v^g_1 + v^g_2 + \cdots + v^g_{2g}) = w_{2g+1}
\end{equation*}
\end{proof}

\subsection{Existence of $2g+1$ pairwise intersecting curves}
We can explicitly describe a $1$-system with $2g+1$ curves all intersecting pairwise.
\begin{theorem} \theolab{canonical-curves}
For any genus $g \geq 1$ closed surface, there is a $1$-system of $2g+1$ curves which all pairwise intersect.
\end{theorem}
\begin{proof}
A genus $g$ surface is homeomorphic to a regular $4g$-gon with opposite sides identified.  We obtain $2g$ simple closed curves from
simple arcs connecting opposite sides of the $4g$-gon.  After identifying sides, all vertices become identified, so we can add a diagonal
to the configuration as well.
\end{proof}
Later, we will use the following slightly stronger statement which is clear from the construction.
\begin{corollary} \corlab{canonical-curves-same-point}
The $1$-system from \theoref{canonical-curves} can be chosen so that all curves pairwise intersect at the same point on the surface.
\end{corollary}

In both genus $1$ and $2$, it can be shown that, up to the mapping class group, there is only one maximal
configuration of curves all of which pairwise intersect exactly once.  In genus $2$, the configuration corresponds to a
star graph in the quotient under hyperelliptic involution.  It would be interesting to know
if this generalizes\footnote{While this paper was under review, Tarik Aougab \cite{A12} showed
that the answer is ``no.''}.
\begin{question}
In higher genus surfaces, is there only one mapping class orbit of $1$-systems such that each pair of curves intersects?
\end{question}

\section{Upper bounds for surfaces with boundary}
We now prove the upper bound of \theoref{boundaryasymptotics}.  First, we improve \propref{all-intersecting}
\begin{lemma} \lemlab{all-intersecting-boundary}
The statement of \propref{all-intersecting} holds for a genus $g$ surface with any number of boundary components.
\end{lemma}
\begin{proof}
Suppose $X$ is a $k$-system in a genus $g$ surface with $n$ boundary components with any pair of
curves intersecting an odd number of times.  Glue in discs into all the boundary components to obtain a closed surface $S$.
Since each curve intersects any other curve transversely an odd number of times, no pair of curves of $X$ can be pairwise isotopic
in $S$.  Consequently, $X$ is a $k$-system satisfying the hypothesis of \propref{all-intersecting} and hence $|X| \leq 2g+1$.
\end{proof}

\subsection{Proof of \theoref{boundaryasymptotics}}
We argue by induction and show that
$$N(1, g, n) \leq (2g+1) + N(1, g, n-1)$$

Let $X$ be a maximal $1$-system in $S_{g, n}$, the surface of genus $g$ and $n$ boundary components.  Glue in a disc $D$
to one of the boundary components to obtain $S_{g, n-1}$.  Since intersection numbers did not increase, curves in $X$ still
pairwise intersect at most once, but some curves may have become isotopic as a result of gluing in the disc.  We will show that a
$1$-system $X'$ on $S_{g, n-1}$ may be obtained by removing at most $2g+1$ curves from $X$.  Then, we have
$$|X| \leq |X'| + 2g+1 \leq N(1, g, n-1) + 2g + 1$$

If two curves become isotopic after gluing in $D$, then since $X$ is a $1$-system, they are disjoint and hence bound an annulus
which necessarily contains $D$.  In particular, any curve can become isotopic to at most $1$ other curve, so we must understand
how many pairs of isotopic curves can occur.  If $\gamma_1,\gamma_2$ and $\eta_1, \eta_2$ are two such pairs, then both pairs
bound an annulus containing $D$ and the annuli must intersect.  Consequently, each curve in an isotopic pair intersects every curve
in the other pairs.  If we construct a set $Y$ by taking one curve from each pair, then $Y$ is as in \lemref{all-intersecting-boundary} and
so $|Y| \leq 2g+1$.  Thus $X' = X \setminus Y$ is the $1$-system on $S_{g, n-1}$ as desired.

In the case of the torus, we have $N(1, 1, 1) \leq N(1, 1, 0)$.  Indeed, if two curves on a one-holed torus become isotopic after gluing in the disc $D$,
then, as before, they bound an annulus containing $D$; however, since the surface is a torus, they also bound an annulus, not containing
$D$, on the other side, and so they must have already been isotopic in the one-holed torus.
\eop

In the next section, we will see that obtaining our lower bounds essentially amounts to reversing the argument in the previous paragraph.

\section{Lower bounds}\seclab{lower}
In this section we prove the lower bounds of \theoref{generalbounds} and \theoref{boundaryasymptotics}. We will first prove \theoref{boundaryasymptotics}.  The lower bound of
\theoref{generalbounds} then follows easily by attaching handles.

\subsection{Proof of lower bound in \theoref{boundaryasymptotics}}
Fix $g \geq 2$.  We show by induction on $n$ that there is a $1$-system $X$ on $S_{g, n}$, the closed surface of genus $g$ with $n$ boundary components, such that
\begin{itemize}
\item $|X| = (2g+1) (n+1)$
\item There is a subset of $Y \subset X$ of size $2g+1$ curves such that those curves and only those curves pairwise intersect all at the same point $\ast$ on the surface.
\end{itemize}
The base case is \corref{canonical-curves-same-point}.  Suppose we have constructed $X$ for $S_{g, n}$.  Let
$D$ be a small disc which contains $\ast$ and intersects only the curves from $Y$.  Via a homeomorphism, we can identify
$D$ with the standard unit disc in $\R^2$ and the arcs from $Y$ as straight-line diagonals all intersecting at the center $= \ast$, none of which is vertical.  Remove
a small disc $D'$ directly above $\ast$, and construct $2g+1$ new simple closed curves as follows.

\begin{figure}[htbp]
\begin{center}
\includegraphics[width=0.5\textwidth]{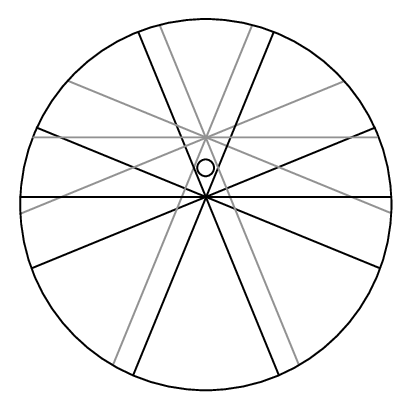}
\end{center}
\caption{The disc $D$ with arcs from $Y$ and their vertical translates (indicated in gray) separated by the disc $D'$.}
\label{figure:lowerbound}
\end{figure}

For each arc $a$ in $D$, place a new arc $a'$ in $D$ parallel to $a$ but above $D'$.  See figure \ref{figure:lowerbound}.
Obtain a simple closed curve $\alpha'$ by continuing $a'$ outside of $D$ along the
curve $\alpha$ in $Y$ which contains $a$.  This can be done over all $\alpha \in Y$ so that:
\begin{itemize}
\item For all $\alpha, \beta \in Y$, the corresponding new curves $\alpha', \beta'$ deposit straight lines in
$D$ and otherwise outside of $D$, the curves $\alpha, \beta, \alpha', \beta'$ are pairwise disjoint.  Hence,
pairwise intersections among $\alpha, \beta, \alpha', \beta'$ are at most $1$.
\item For all $\alpha \in Y$ and $\beta \in X \setminus Y$, the new curve $\alpha'$ has intersection
$i(\alpha', \beta) = i(\alpha, \beta)$
\end{itemize}
Let $X'$ be $X$ with the newly constructed curves added.  Each new curve $\alpha'$ bounds with its
``parent'' $\alpha$ an annulus containing the (removed) disc $D'$.
Since $g \geq 2$, the complement of the annulus has genus at least $1$, and thus is not an annulus; consequently
$\alpha$ and $\alpha'$ are non isotopic.  Furthermore, $\alpha'$ is not isotopic to any other curve in $X'$ since
if it were, then, after replacing $D'$, we would see that $\alpha$ were isotopic to some other curve in $X$,
a contradiction.  Note that $X'$ has exactly $2g+1$ more curves and $Y\subset X'$ still has the desired properties.

Notice that this argument fails in $g = 1$ only because if $n = 0$, then the complement of the annulus containing
$D$ would in fact be another annulus. However, once $n \geq 1$, the complement would contain a boundary component,
and the argument proceeds mutatis mutandis, but with the smaller lower bound of $(2g+1) n$.
\eop

\subsection{Proof of the lower bound in \theoref{generalbounds}}
The cases $g = 2, 3$ were established in \secref{gen2} so assume $g \geq 4$.
Let $m = g/2 \geq 2$ if $g$ is even and $m = (g-1)/2$ if $g$ is odd and let $n = g-m$.  By \theoref{boundaryasymptotics},
there is a $1$-system $X$ of size $(2m+1)(2n+1)$ on a surface of genus $m$ and $2n$ boundary components.  Gluing $n$ handles
to the $2n$ boundary components does not cause curves to become isotopic, the system consisting of $X$ and the $n$ curves
going around the handles is a $1$-system of at least $g^2 + \frac{5}{2}g$ curves on a genus $g$ surface.
\eop

\bibliographystyle{spmpscinat}

\end{document}